\documentclass[10pt,refreqno,oneside]{amsart}
\usepackage{amssymb}

\newtheorem{theorem}{Theorem}
\newtheorem{proposition}[theorem]{Proposition}

\newtheorem{lemma}[theorem]{Lemma}

\newcommand{\cnums}{\mathbb{C}}
\newcommand{\reals}{\mathbb{R}}

\newcommand{\bu}{\mathbf{u}}
\newcommand{\bx}{\mathbf{x}}

\newcommand{\bI}{\mathbf{I}}
\newcommand{\bxi}{\boldsymbol{\xi}}

\newcommand{\rpx}{\reals[\bx]}
\newcommand{\rpxi}{\reals[\bxi]}

\newcommand{\cpxi}{\cnums[\bxi]}
\newcommand{\cpx}{\cnums[\bx]}
\newcommand{\cpI}{\cnums[\bI]}
\newcommand{\cfx}{\cnums(\bx)}
\newcommand{\cfxi}{\cnums(\bxi)}
\newcommand{\cfI}{\cnums(\bI)}

\newcommand\lb{\left[}
\newcommand\rb{\right]}
\newcommand\lp{\left(}
\newcommand\rp{\right)}

\title[Reflection quotients]{Reflection quotients in Riemannian
  Geometry. A Geometric Converse to Chevalley's Theorem} 

\author{R. Milson}


\address{Dept of Math., Dalhousie U., Halifax, Canada, B3H 3J5}

\email{milson@mscs.dal.ca} 

\thanks{Supported by NSERC grant 228057.}

\begin{document}
\begin{abstract}
  Chevalley's theorem and it's converse, the Sheppard-Todd theorem,
  assert that finite reflection groups are distinguished by the fact
  that the ring of invariant polynomials is freely generated.  We show
  that in the Euclidean case, a weaker condition suffices to
  characterize finite reflection groups, namely that a
  freely-generated polynomial subring is closed with respect to the
  gradient product.
\end{abstract}
\maketitle

\pagestyle{myheadings}

From the standpoint of invariant theory, finite reflection groups are
distinguished by the property that the corresponding algebra of
invariants is freely generated.  The existence of free generators of
the invariant algebra is known as Chevalley's theorem
\cite{chevalley}.  The converse, i.e.  the statement that a finite
group with a freely generated invariant algebra is necessarily
generated by reflections, is known as the Sheppard-Todd theorem
\cite{sheptodd}.  The customary proofs of these results use algebraic
methods.

The purpose of this note is to propose an alternate characterization
of finite reflection groups over the reals, and to prove the result
using the language and ideas of Riemannian geometry.  Riemannian
theory is relevant because a finite reflection group over $\reals$
determines a Euclidean structure (this can be seem by a simple
averaging argument), and therefore, without loss of generality, the
elements of the reflection group can be assumed to be Euclidean
automorphisms.  Now the structure of Riemannian geometry is
specified by a fundamental covariant: the gradient
operation.  Axiomatically then, the gradient operation is preserved by
all Riemannian automorphisms. 

Specializing to Euclidean space, if $P(x^1,\ldots,x^n)$,
$Q(x^1,\ldots, x^n)$ are two polynomials that are invariant with
respect to some Euclidean reflections, then the corresponding gradient
product
$$\nabla P\cdot \nabla Q = \sum_{i} \frac{\partial
  P}{\partial x^i} 
\frac{\partial Q}{\partial x^i}$$
is also an invariant polynomial.  It
follows immediately that, if $I^1, \ldots, I^n$ is a free basis of the
invariant algebra, then the matrix of gradient cross products also
consists of polynomials in the basis elements, i.e.
\begin{equation}
  \label{eq:gradcprod}
\nabla I^i \cdot \nabla I^j = g^{ij}(I^1,\ldots,I^n),\quad
i,j=1,\ldots,n   
\end{equation}
where the $g^{ij}$ are $n$-variable polynomials.  

We will show that this property characterizes the class of finite
reflection groups in the following sense. 
\begin{theorem}
\label{thrm:main}
Let $I^1, \ldots, I^n$ be algebraically independent, homogeneous, real
polynomials in $n$ variables, and let $G$ be the group of linear
automorphisms that leaves them invariant.  If all the corresponding
Euclidean gradient cross-products are themselves polynomial in
$I^1,\ldots,I^n$, then $I^1,\ldots,I^n$ freely generate the algebra of
$G$-invariants.
\end{theorem}
\noindent
Accepting this theorem as true, the Shepard-Todd theorem then
implies that $G$ is generated by reflections.

The restriction of the field to $\reals$ and the assumption of
Euclidean signature are indespensible;  the theorem is not true
without them.  Consider the following 2-dimensional counterexamples
$$ I^1 = x^1+ix^2,\quad I^2 = (x^1)^2+(x^2)^2.$$
The matrix of gradient cross-products is:
$$\lp\begin{array}{ll}
  0 & 2I^1 \\ 2I^1 & 4I^2
\end{array}\rp
$$
However, $x^1-ix^2 = I^2/I^1$ and hence $G$
is trivial.  The conclusion of the theorem does not hold.
Alternatively, one can change the above into a real counterexample
based on non-Euclidean signature.  Take
$$I^1 = x^1+x^2,\quad I^2 = (x^1)^2-(x^2)^2.$$

The key to the proof of Theorem \ref{thrm:main} is to regard the
matrix of gradient products, $g^{ij}$, as the contravariant form of a
Riemannian metric tensor.  The matrix of gradients cross-products has
been studied by in the context of singular projections
\cite{arnold}\cite{saito}.  To the author's best knowledge,
the Riemannian-geometric interpretation of $g^{ij}$ is novel.

There is a good reason why others may have hesitated to interpret the
matrix of gradient cross-products as a Riemannian tensor.  The idea
cannot be made entirely rigorous because of an essential complication.
Let $\delta$ denote the discriminant
$$\delta = \det g^{ij}, $$ and let us note that
the tensor $g^{ij}$ is degenerate at points
where $\delta=0$. The presence of degeneracies means that the
$g^{ij}$ do not define a Riemannian structure in the usual sense of
this term.

To put it another way, the $g^{ij}$ permit us to regard the map
$$
\Pi: \reals^n \rightarrow \reals^n, $$
where $\Pi =
(I^1,\ldots,I^n)$, as an isometry. However, there is a complication:
$\Pi$ is not a regular covering.  Note that $\Pi^*\delta$ is the
square of the Jacobian of $\Pi$,
$$J = \det \lb \frac{\partial I^j}{\partial x^i}\rb,$$
and hence the equation $\delta =0$ picks out the images of points where
this map has less than maximal rank.  

The preceding remarks should be taken as an indication that the
classical Riemannian theory admits an interesting generalization.  The
generalized theory admits objects that are quotients of Riemannian
manifolds by reflection groups --- a kind of Riemannian
reflection-orbifold if you will.  However, the local covering data for
such orbifolds does not have to be specified as a separate item of
information; it is encoded implicitly by the degeneracies of the
metric tensor $g^{ij}$.  The maps of the theory are more general than
the usual immersions, isometries and submersions of Riemannian
geometry, and include variable rank ``foldings''.

We will relegate the development of such a theory to subsequent
publications.  Some preliminary remarks are available in \cite{milson}.
Also, for a different type of geometric theory based on degenerate
$g^{ij}$ see \cite{hermann}.  In this note we restrict our attention to the
proof of the Theorem \ref{thrm:main}, and merely use the geometric
ideas as a guiding principle.  The fact of the matter is that the
restriction of the $g^{ij}$ to the domain $\{\delta\neq0\}$ defines a
perfectly standard, albeit incomplete, Riemannian structure.  This
observation underlies and motivates the proof technique.

We now turn to the proof of our theorem.  Let us fix notation and
hypotheses.  Let $I^1, \ldots, I^n$ be algebraically independent, real
homogeneous polynomials in $n$ indeterminates $x^1,\ldots, x^n$.  It
will be convenient to carry out the proof over the field of complex
numbers, and so we identify the polynomials in question with maps from
$\cnums$ to $\cnums$, which commute with complex conjugation. In the
same spirit, we regard $\Pi = (I^1,\ldots, I^n)$ as a map from
$\cnums^n$ to $\cnums^n$ and define $G$ to be the group of complex
linear automorphisms that preserve $\Pi$.

We shall write $\cpx$ and $\cpI$ for the polynomial algebras
generated, respectively, by $x^1,\ldots,x^n$, and by $I^1,\ldots,I^n$.
We shall write $\cfx$ and $\cfI$ for the corresponding fraction
fields.  It will also be convenient to introduce $n$ additional
indeterminates $\xi^1,\ldots, \xi^n$ to serve as the coordinates of
the codomain:
$\Pi^*(\xi^j) = I^j.$
Throughout we assume that the $\cpI$
is closed with respect to the gradient product, i.e. that there exist
polynomials $g^{ij}\in\rpxi$ such that \eqref{eq:gradcprod} holds.  We
can now define a  gradient operation on $\cpxi$:
$$  \nabla \alpha \cdot \nabla \beta = \sum_{ij} \frac{\partial
  \alpha}{\partial \xi^i} 
\frac{\partial \beta}{\partial \xi^j}\, g^{ij},\quad \alpha,\beta\in\cpxi.
$$
This operation is,
by construction, compatible with the gradient operation on $\cpx$:
\begin{equation}
  \label{eq:gradcompat}
\Pi^*(\nabla \alpha \cdot \nabla \beta) = \nabla (\Pi^*\alpha) \cdot
\nabla(\Pi^*\beta),
\end{equation}
for all $\alpha, \beta\in\cpxi$.
We
let $\delta\in\rpxi$ and $J\in\rpx$ denote, respectively, the
discriminant and the Jacobian as per above, and note that
$$ \Pi^*\delta = J^2.$$

\begin{proposition}
\label{prop:galois}
The group $G$ of linear $\Pi$-automorphisms is a subgroup of
$\mathrm{O}_n(\cnums)$, and is isomorphic to the group of field
automorphisms of $\cfx$ over $\cfI$.  Furthermore $\cfx$ is a normal
extension of $\cfI$, and hence $G$ is isomorphic to the Galois group
of $\cfx$ over $\cfI$.
\end{proposition}
\begin{proof}
  By definition, each $g\in G$ defines an automorphism of $\cpx$ over
  $\cpI$ and hence of $\cfx$ over $\cfI$.  We therefore have a natural
  inclusion 
  $$G\subset\mathrm{Aut}(\cfx/\cfI).$$
  We now prove that
  this inclusion is, in fact, an isomorphism, as well as show that the
  elements of $G$ are orthogonal.
  
  Let $d$ be the degree of the extension.  Using a primitive element,
  one can straightforwardly show that there exists a dense open subset
  $U\subset \cnums^n$ such that for every $\bxi\in U$ the set
  $\Pi^{-1}(\bxi)$ consists of $d$ distinct preimages.  We avoid the
  locus $\{\delta=0\}$ and choose an open $U_0\subset U$ sufficiently
  small so that $\Pi^{-1}(U)$ is the union of $d$ disjoint open sets
  $V_1,\ldots, V_d$ and such that the restriction of $\Pi$ to each of
  these is non-singular.  We therefore have $d$ complex-analytic maps
  $$\sigma_i: U_0 \rightarrow V_i,\quad i=1,\ldots, d$$
  that are local inverses of $\Pi$.  For each $i=1,\ldots,d$ we define
  $g_i: V_1\rightarrow V_i$ by
  $$g_i = \sigma_i \circ \Pi\Big|_{V_1}.$$
  We then note that each
  $g_i$ preserves the gradient product, and hence defines a local
  automorphism of the dot product structure on the tangent bundle of
  $\cnums^n$.  Such a local automorphism preserves straight lines, and
  hence must be the restriction of a global, semi-linear
  transformation:
  $$g_i\in\mathrm{O}_n(\cnums)\ltimes \cnums^n.$$
  However, the components of $\Pi$ are homogeneous polynomials and
  $$\Pi\circ g_i = \Pi,$$
  by definition.  Hence, each $g_i$ maps the
  origin to itself, hence is linear, and hence is an element of
  $G$.  Thus, we have shown that $G$ contains $d$ elements of
  $\mathrm{O}_n(\cnums)$.  However, the order of the automorphism
  group   cannot exceed $d$, and the desired conclusions follow.
  
  We now prove that $\cfx$ is Galois over $\cfI$ by showing that
  $\cfI$ is the fixed field of $G$.  Let $p\in\cfx$ such that
  $p\notin\cfI$ be given.  The minimal polynomial of $p$ over $\cfI$
  has degree $2$ or more, and hence there will exist points $\bx_1,
  \bx_2\in \cnums^n$ such that $\Pi(\bx_1)=\Pi(\bx_2)$ but such that
  $p(\bx_1)\neq p(\bx_2)$.  However, using the technique in the
  preceding paragraph one can find a $g\in G$ such that $g(\bx_1) =
  \bx_2$, and hence $p$ is not $G$-invariant.  Therefore the fixed
  field of $G$ coincides with $\cfI$.
\end{proof}


 \begin{lemma}
   \label{lemma:onto}
 The open $\{\delta\neq 0\}\subset\cnums^n$ is in the
 image of $\Pi$.
 \end{lemma}
 \begin{proof}
   Note that the open subset in question is path connected.
   Consequently if $\Gamma\subset \{\delta\neq 0\}$ is a continuous
   path between points $\bxi_1, \bxi_2$, and if
   $\bxi_1$ is in the image of $\Pi$, then so is $\bxi_2$.  This is
   because $\Pi$ has a local inverse wherever $J\neq 0$, and hence
   $\Gamma$ can be lifted to a continuous path.  However, at least one
   such $\bxi_1$ is bound to exist, and the desired conclusion follows.
 \end{proof}

\begin{proposition}
  \label{prop:trapping2}
  If $\lambda\in\cpxi$ is such that $\delta=0$ wherever $\lambda=0$,
  then $\nabla \log \lambda$ is a well defined derivation of $\cpxi$.
  In other words for every $\rho\in\cpxi$ there exists a
  $\sigma\in\cpxi$ such that
  $$\nabla \lambda(\rho) = \sigma \lambda.$$
\end{proposition}
\begin{proof}
  It suffices to give a proof for irreducible $\lambda$.  Suppose the
  proposition is false.  Then, there exists a $\rho\in\cpxi$ such that
  $\nabla \rho$ is transverse to the subvariety $\{\lambda=0\}$.  In
  other words, $\nabla \rho(\lambda)$ does not vanish identically on
  $\{\lambda=0\}$.  The covariant formula for the Laplacian gives
  \begin{equation}
    \label{eq:lapform}
    \Delta \rho = \sum_i \frac{\partial}{\partial\xi^i} (\nabla
      \xi^i(\rho))  - \frac12 (\nabla \log \delta)(\rho).
  \end{equation}
  It follows that
  $$\Pi^*\lp(\nabla\log\delta)(\rho)\rp\in\cpx,$$ 
  and since $\Pi^*\lambda$ divides $\Pi^*\delta$
  $$\Pi^*\lp(\nabla\log\lambda)(\rho)\rp\in\cpx$$
  as well.

  Next, let $\Phi_t$ be the
  one-parameter flow generated by $\nabla \rho$, and $\phi_t$ the one
  parameter flow generated by $\nabla (\Pi^*\rho)$.  By
  \eqref{eq:gradcompat}, the two flows intertwine with $\Pi$:
  \begin{equation}
    \label{eq:pigradflow}
    \Phi_t\circ\Pi = \Pi\circ\phi_t.
  \end{equation}
  Choose a $\xi_0\in\{\lambda=0\}$ such that
  $$\nabla\rho(\lambda)_{\bxi_0}\neq 0.$$
  It follows immediately that
  $$\nabla\rho(\delta)_{\bxi_0}\neq 0,$$
  and hence, for all sufficiently
  small $t$ we have 
  $$\delta(\Phi_t(\xi_0))\neq 0.$$
  By the Lemma,
  $\Phi_t(\xi_0)$ is in the image of $\Pi$ , and hence by
  \eqref{eq:pigradflow} so is $\xi_0$.  However 
  $$\Pi^*(\nabla\rho(\lambda)) = \Pi^*(\lambda)
  \Pi^*\lp(\nabla\log\lambda)(\rho)\rp.$$ 
  The right hand side is zero at the pre-image of $\xi_0$,
  and hence so it the left hand side --- a contradiction.
\end{proof}


\begin{proposition}
  \label{prop:orthogirred}
If $\lambda_1,\lambda_2\in\cpxi$ are distinct irreducible factors of
the discriminant $\delta$, then 
$$\nabla\lambda_1\cdot \nabla\lambda_2= 0.$$
\end{proposition}
\begin{proof}
  By the preceding proposition, $\nabla\lambda_1\cdot
  \nabla\lambda_2$ is divisible by $\lambda_1 \lambda_2$. However, the
  degree of $\Pi^*$ of the former is smaller than the degree of
  $\Pi^*$ of the latter, and the desired conclusion follows.
\end{proof}
\begin{proposition}
\label{prop:detisharmonic}
  The Jacobian $J$ is a harmonic polynomial.
\end{proposition}
\begin{proof}
  The condition $\Delta J = 0$ is equivalent to 
  \begin{equation}
    \label{eq:deltalogdelta}
    \Delta \log\delta = - \frac{(\nabla \delta)(\log\delta)}{2\delta}.   
  \end{equation}
  We show that the latter equation holds.  Write 
  $$\nabla \log\delta = \sum_i \alpha_i \frac{\partial}{\partial
    \xi^i},$$
  where 
  $$\alpha_i = \nabla \xi^i(\log \delta)=(\nabla
  \log\delta)(\xi^i).$$
  By the preceding Proposition the $\alpha_i\in\cpxi$, and since
    degree of $\Pi^*\alpha_i$ must be lower than the degree of
    $\Pi^*\xi^i$ we must have
    $$\frac{\partial\alpha_i}{\partial \xi^i} = 0$$
    for all $i$.   Equation \eqref{eq:deltalogdelta} now follows from
    the formula for the Laplacian shown in \eqref{eq:lapform}.
\end{proof}


\begin{proposition}
\label{prop:intext}
  $\cpx$ is an integral extension of $\cpI$.
\end{proposition}
\begin{proof}
  Suppose not.  Consider a $p\in\cpx$ that is not integral over
  $\cpI$.  The coefficients of the minimal polynomial over $\cfI$ are
  symmetric functions of the $G$-conjugates of $p$, and hence, by
  Proposition \ref{prop:galois} are polynomials in $\cpx$.  We have
  now demonstrated that if the extension is non-integral, then there
  must exist non-integral elements of degree 1, i.e. $r\in\cpx$ such
  that 
  $$r=\Pi^*(\alpha/\beta)$$
  for some relatively prime $\alpha,
  \beta\in\cpxi$. Choose one such $r$. The repeated application of the
  Laplacian will eventually reduce every polynomial to a constant, and
  hence without loss of generality we may assume that $\Delta r$ is a
  constant.  We then have the following identity in $\cfxi$:
  \begin{equation}
    \label{eq:intext1}
    \beta\lp\Delta\alpha - 2(\nabla \log\beta)(\alpha) - \beta\Delta
    r\rp = \alpha ( \Delta \beta - 2(\nabla \log\beta)(\beta)).
  \end{equation}
  We can also assert that $\delta=0$ wherever $\beta=0$.  If not,
  there would exist a $\bxi_0$ such that
  $$\alpha(\bxi_0)\neq 0,\quad \beta(\bxi_0) =0,\quad
  \delta(\bxi_0)\neq 0.$$
  However, then by Lemma \ref{lemma:onto},
  $\bxi_0=\Pi(\bx_0)$ for some $\bx_0\in\cnums^n$, and hence
  $$\alpha(\bxi_0) = \beta(\bxi_0) r(\bx_0) = 0,$$
  a contradiction.
  
  Hence, by Proposition \ref{prop:trapping2}, $\nabla \log\beta$ is a
  well-defined derivation of $\cpxi$.  It follows that
  \eqref{eq:intext1} is, in fact, a relation between polynomials.
  Since $\alpha$ and $\beta$ are relatively prime, equation
  \eqref{eq:intext1} implies that $\beta$ divides
  $$\Delta \beta - 2(\nabla \log\beta)(\beta).$$
  However, the degree of $\Pi^*$ of this expression is smaller than
  the degree of $\Pi^*\beta$, and hence
  \begin{equation}
    \label{eq:intext2}
    \Delta \beta = 2(\nabla \log\beta)(\beta).   
  \end{equation}
  
  Next, let $\lambda$ be an irreducible factor of $\beta$; let's say
  with multiplicity $j$.  We have shown that $\lambda$ is also an
  irreducible factor of $\delta$; let's say with multiplicity $k$.  By
  Proposition \ref{prop:orthogirred} all irreducible factors of
  $\delta$ are mutually orthogonal, and hence equation
  \eqref{eq:intext2} implies that
  \begin{equation}
    \label{eq:intext3}
    \Delta \lambda = (j+1) (\nabla \log\lambda)(\lambda).
  \end{equation}
  However in Proposition \ref{prop:detisharmonic} we showed that
  $$  \Delta \delta = \frac12 (\nabla \log\delta)(\delta),$$
  and hence
  \begin{equation}
    \label{eq:intext4}
    \Delta \lambda = \lp 1-\frac{k}2\rp (\nabla \log\lambda)(\lambda).    
  \end{equation}
  Combining \eqref{eq:intext3} and \eqref{eq:intext4} we conclude that
  $$\Delta \lambda = \Delta \lambda^2 = 0.$$
  
  There does not exist a real polynomial that obeys such relations.
  However, $\delta$ is a polynomial with real coefficients, and hence
  $\bar{\lambda}$ is also an irreducible factor of $\lambda$ that
  obeys the same relations.  By Proposition \ref{prop:orthogirred},
  $\lambda$ and $\bar{\lambda}$ are mutually orthogonal, and hence
  $$\Delta \lp\lambda \bar{\lambda} \rp = \Delta \lp \lambda \bar\lambda
  \rp^2 = 0.$$  This is impossible, and the desired conclusion follows.
\end{proof}
We are now ready to give the proof of Theorem \ref{thrm:main}.
Suppose $r\in\cpx$ is a $G$-invariant.  Since $G$ is isomorphic to the
Galois group of $\cfx$ over $\cfI$, there exists a $\rho\in\cfxi$ such
that $r=\Pi^*\rho$.  However, we have just shown that the extension is
integral, and so in fact $\rho\in\cpxi$.   This proves the Theorem.

\end{document}